\documentclass[11pt]{amsart}

\usepackage{color}
\usepackage{graphicx}
\usepackage{color}
\usepackage[latin1]{inputenc}
\usepackage{psfrag}

\usepackage[OT1]{fontenc}
\usepackage{pslatex}

\usepackage{a4wide}
\usepackage{amsfonts}
\usepackage{amsmath}
\usepackage{amssymb}
\usepackage{amsthm}	
\usepackage{enumerate}

\usepackage{bbm}
\theoremstyle{remark}
\newtheorem{rem}{Remark}[section]

\numberwithin{equation}{section}

\theoremstyle{definition}
\newtheorem{teo}{Theorem}[section]

\newtheorem{lema}[teo]{Lemma}
\newtheorem{prop}[teo]{Proposition}
\newtheorem{cor}[teo]{Corollary}

\newtheorem{defi}[teo]{Definition}

\newcommand{\R}{{\mathbb R}}

\newcommand{\N}{{\mathbb N}}

\renewcommand{\d}{{\rm d}}

\newcommand{\ve}{{{\varepsilon}}}

\def\R{\mathbb{R}}

\def\W{\mathcal W}
\def\F{\mathcal F}

\def\1{{\bf 1}}

\begin{document}

\title{Small random perturbations of a dynamical system with blow-up}

\author{Pablo Groisman and Santiago Saglietti}

\thanks{ 2000 {\it Mathematics Subject Classification:} 60H10, 34C11.}

\keywords{Random perturbations, explosions, stochastic differential equations, blow-up.}

\address{Departamento de Matem\'atica\hfill\break
\indent Facultad de Ciencias Exactas y Naturales\hfill\break
\indent Universidad de Buenos Aires \hfill\break
\indent Pabell\'{o}n I, Ciudad Universitaria \hfill\break
\indent C1428EGA Buenos Aires, Argentina.
}
\email{{\tt pgroisma@dm.uba.ar, ssaglie@dm.uba.ar}
\hfill\break\indent {\it
Web page:} {\tt http://mate.dm.uba.ar/$\sim$pgroisma}}

\date{}

\begin{abstract}
We study small random perturbations by additive white-noise of a spatial
discretization of a reaction-diffusion equation with a stable
equilibrium and solutions that blow up in finite time. We prove that the
perturbed system blows up with total probability and establish its order of
magnitude and asymptotic distribution. For initial data in the domain of
explosion we prove that the explosion time converges to the deterministic one
while for initial data in the domain of attraction of the stable equilibrium we
show that the system exhibits metastable behavior.
\end{abstract}

\maketitle

\section{Introduction}
We consider small random perturbations of the following ODE
\begin{equation}
\label{1.1}
\left\{\begin{array}{lcll}
U'_1 &= &\frac{2}{h^2} ( -U_1 + U_2 ),\\
U'_i &= &\frac{1}{h^2} ( U_{i+1} - 2U_i + U_{i-1} ) &\,\,\, 2 \leq i \leq d-1,\\
U'_d &= &\frac{2}{h^2} ( -U_d + U_{d-1} +hg(U_d) ).
\end{array}\right.
\end{equation}
Here $g\colon \R \to \R$ is a reaction term given by $g(x) = (x^+)^p -x$ with $p>1$, and $h>0$ is a parameter. We also impose an initial condition $U_0 \in \R^d$. This kind of systems arise as spatial discretizations of diffusion equations with nonlinear boundary conditions of Neumann type. In fact, it is known that as $h\to 0$ solutions to this system converge to solutions of the PDE

$$
\left \{\begin{array}{rcll}
u_t(t,x) & = & u_{xx}(t,x) & 0<x<1, 0\le t<T,\\
u_x(0,t) & = & 0 & 0\le t <T,\\
u_x(1,t) &= & g(u(1,t)) & 0\le t <T,\\
u(x,0) & = & u_0(x) & 0\le x \le 1.
\end{array} \right.
$$

This and more general reaction-diffusion problems including for instance the
possibility of a nonlinear source term like $g$ and other type of boundary
conditions appear in several branches of pure and applied mathematics. They have
been used to model heat transfer, exothermic chemical reactions, population
growth models, geometric flows, etc.

An important feature of this type of problems is that they
admit solutions which are local in time, with the possibility of blow-up in
finite time. The asymptotic behavior of solutions to \eqref{1.1} can be briefly
summarized as follows (we give a detailed description afterwards): the system
has two equilibriums $U_0\equiv 0$ and $U_0 \equiv 1$. The first one is stable
while the second is unstable. Hence, there exists a domain of attraction $D_0$
for the zero solution such that if $U_0 \in D_0$ then the solution
$U(t)=(U_1(t), \dots, U_d(t))$ with initial condition $U_0$ is globally defined
and $U(t) \to 0$ as $t \to \infty$. There exists also a stable manifold for the
unstable equilibrium which is of co-dimension one and coincides with the
boundary of $D_0$. For $U_0 \in {\overline D_0}^c$ the solution $U$ blows up in
finite time $T=T(U_0)$.

Since mathematical models are not exact, it is important to understand what
changes arise in the behavior of the system when it is subject to
perturbations. We study random perturbations given by additive white-noise. More
precisely, we consider Stochastic Differential Equations (SDE) of the form
\begin{equation}
\label{1.epsilon}
\left\{\begin{array}{rcll}
\d U^{u,\ve}_1 &= &\frac{2}{h^2} (-U^{u,\ve}_1 + U^{u,\ve}_2 )\,\d t + \ve \d W_1,\\
\d U^{u,\ve}_i &= &\frac{1}{h^2} (U^{u,\ve}_{i+1} - 2U^{u,\ve}_i + U^{u,\ve}_{i-1})\, \d t + \ve {\rm d}W_i &  2 \leq i \leq d-1,\\
\d U^{u,\ve}_d &= &\frac{2}{h^2} ( -U^{u,\ve}_d + U^{u,\ve}_{d-1}
+hg(U^{u,\ve}_d) )\, \d t + \ve \d W_d,& \\
\end{array}\right.
\end{equation}
which can be written in matrix form as
\begin{equation}
\label{A}
dU^{u,\ve} = (-AU^{u,\ve} + \frac2h g(U_d^{u,\ve})e_d)\,dt + \ve dW.
\end{equation}
Here $W=(W_1, \dots W_d)$ is a $d-$dimensional standard Brownian motion,
$\ve>0$ is a small parameter and $e_d=(0,\dots,1)$ is the $d$-th canonical
vector on $\R^d$. In the sequel we use $U^{u,\ve}$ for a solution to
\eqref{1.epsilon} with initial condition $U^{u,\ve}(0)=u \in \R^d$. In the case
$\ve=0$ we are left with the deterministic equation and so we use the notation
$U^u:=U^{u,0}$ to denote a solution to \eqref{1.1}.

The field $b(U):= -AU + \frac2hg(U_d)e_d$ is a gradient ($b=-\nabla \phi$) with
potential given by
$$
\phi(U) = \frac{1}{2} \langle AU , U \rangle - \frac{2}{h}\Big(\frac{\;\;\;|U_d^+|^{p+1}}{p+1} - \frac{\,\,{U_d}^2}{2}\Big).
$$
The SDE associated to this energy functional can be compared with the classic
double-well potential model, which we now briefly summarize. We refer to
\cite[p. 294]{OV} for a more detailed description.

In the double-well potential model one considers a stochastic differential equation of the form
\begin{equation}
\label{eq.general}
dX^{\varepsilon} = r(X^{\varepsilon})\,dt + \varepsilon\,dW
\end{equation}
where $W$ is a standard $d$-dimensional Brownian motion and $r$ is a globally
Lipschitz gradient field over $\R^d$ given by the double-well potential
$\tilde{\phi}$. More precisely, this
potential $\tilde\phi$ possesses exactly three critical points: two local
minima $p$ and $q$ of different depth and a saddle point $z$ with
higher energy, that is $\tilde\phi(z) > \tilde\phi(p)>\tilde\phi(q)$ .
Each minimum corresponds to a stable equilibrium and hence for initial data
lying outside the stable manifold of $z$, the deterministic system ($\ve =0$)
converges to one of them depending on the initial condition. When considering
random perturbations, for compact time intervals the stochastic
system converges as $\ve \to 0$ to the deterministic one uniformly but the
qualitative
behavior of the perturbed system is quite different from that of the
deterministic solution for large times. If the potential grows
fast enough at infinity the resulting stochastic system admits a stationary
probability measure which converges to a Dirac delta concentrated
at the bottom of the deepest well $q$. Hence, for initial data in the domain of
attraction of the shortest well $p$ we observe that
\begin{enumerate}
\item [(i)] Due to the action of the field $r$, the process is attracted
towards the bottom of the shortest well $p$; once near $p$, the field
becomes negligible and the process is then pushed away from the bottom of the
well by noise. Being apart from $p$, noise becomes overpowered by the field
$r$ and this allows for the previous pattern to repeat itself: a large number of
attempts to escape from the given well, followed by a strong attraction
towards its bottom. This phase is known as {\em thermalization}.
\item [(ii)] Eventually, after many frustrated attempts, the process succeeds
in overcoming the barrier of potential and reaches the deepest well.
Since the probability of such an event is small, we expect this {\em tunneling
time} to be exponentially large. Moreover, due to the large number of attempts
that are necessary, we expect this time to show little memory.
\item [(iii)] Once in the deepest well, the process behaves as in (i). Since
the new barrier of potential is higher, the next tunneling time is expected to
happen on a larger time scale.
\end{enumerate}
This description was proved rigorously in \cite{Day, GOV, MOS, FW, MS} using
different techniques. The phenomenon is known as {\em metastability}. For a
detailed description of it we refer to \cite{OV}.

Coming back to our potential $\phi$, the situation is slightly more
complex. Instead of having a deepest well, we have a direction along which the
potential goes to $-\infty$ and, hence, the size of the ``deepest well'' is now
infinity and there is no return from there. Moreover, since the potential
behaves like $-s^{p+1}$ in this direction, if the system falls in this ``well'',
it reaches infinity in finite time (explosion).

The purpose of this paper is to study the metastability phenomenon for this
kind of potentials where there is a shortest (finite) well and a deepest well
which yields to infinity in finite time. The ideas developed here can be
extended to other systems with the same structure. The typical situation with
this kind of geometry is the case of reaction-diffusion equations where the
reaction comes from a nonlinear source with superlinear behavior at infinity
such as
$$
u_t = u_{xx} + u_+^p,
$$
with $p>1$, in a bounded domain of $\R$ and homogeneous Dirichlet boundary
conditions. In this case the diffusive term pushes the solution towards zero (a
stable equilibrium) while the source $u_+^p$ pushes it to infinity. In this
situation we expect the same behavior as the one of solutions to
\eqref{1.epsilon}.

Since the drift in \eqref{1.epsilon} is not globally Lipschitz, we are only able
to prove the existence of local solutions and in fact, explosions occur for
solutions of \eqref{1.epsilon}. In particular, classical large deviation
principles as well as other Freidlin-Wentzell estimates do not apply directly.
All of these results deal with globally Lipschitz coefficients. Also, the loss
of memory for the tunneling time was proved only in the globally Lipschitz case
where explosions do not occur. The only exception is the work of
Azencott \cite{A} where locally Lipschitz coefficients are considered and
explosions are allowed, but the large deviations estimates developed there apply
only to neighborhoods of solutions which do not explode in a fixed time interval
(and hence the perturbed system is automatically defined in the whole interval
for $\ve$ small enough). In that work the author also considers the exit from a
domain problem, but explosions are not allowed in  his analysis.

As opposed to this last case, we specifically focus on trajectories that blow
up in finite time. The asymptotic behavior (as $\ve \to 0$) of the explosion
time for \eqref{1.epsilon} is not understood yet, and this is the goal of this
article.

In order to study this kind of systems, localization techniques may be
applied but this has to be done carefully. The main difficulties lie in (i) the
geometry of the potential (and its respective truncations) which is far from
being as simple as in the double-well potential and (ii) the explosion phenomena
itself. Localization techniques apply reasonably well to deal with the process
until it escapes any bounded domain, but dealing with
process from there up to the explosion time requires different tools, which
include a careful study of the blow-up phenomenon. Clearly,
localization arguments are useless for this last part. 

The paper is organized as follows. In Section \ref{defr} we give the
necessary definitions, review some Freidlin-Wentzell estimates and
detail the results of this article. Section \ref{deter} is devoted to giving a
detailed description of the deterministic system \eqref{1.1}. In Section
\ref{estoc.exp} we begin our analysis of the stochastic system. We prove that
explosions occur with probability one for every initial datum. In Section
\ref{estoc.conv} we prove that for initial data in the domain of explosion, the
explosion time converges to the deterministic one as $\ve \to 0$. Finally,
throughout Section \ref{meta} we study the characteristics associated to
metastability for initial datum in the domain of attraction of the origin:
exponential magnitude of the explosion time and asymptotic loss of memory.

\section{Definitions and results}
\label{defr}

\subsection{Solutions up to an explosion time}

\noindent Throughout the paper we study stochastic differential equations of
the form
\begin{equation}\label{eqexplosion0}
dX = \tilde{b}(X)\,dt + \varepsilon\,dW
\end{equation} where $\varepsilon > 0$ and $\tilde{b}:\R^d \longrightarrow
\R^d$ is locally Lipschitz. It is possible that such equations do not admit
strong solutions in the usual sense as these may not be globally defined but
defined \textit{up to an explosion time} instead. We now formalize the idea of
explosion and properly define the concept of solutions for this kind of
equations. We follow \cite{McKean}.

\medskip
\begin{defi}\label{solexp} A \textit{solution up to an explosion time} of the
stochastic differential equation \eqref{eqexplosion0} on the probability space
$(\Omega, \mathcal{F}, P),$ with respect to a filtration $(\F_t)_{t \geq 0}$
satisfying the usual conditions and a fixed Brownian motion $(W_t, \F_t)_{t \geq
0}$ with (a.s. finite) initial condition $\xi$ is an adapted process $X$ with
continuous paths taking values in $\R^d\cup \{\infty\}$ which satisfies the
following properties:

\noindent \begin{itemize}
\item If we define $\tau^n = \inf\{t>0 : |X(t)|=n\}$ then for every $n\ge1$ we have
$$
P\Bigg( \int_{0}^{t\wedge{\tau^n}} |\tilde{b}(X(s))|\,ds < +\infty \Bigg)=1 \hspace{0,5cm} \,\forall\,\,\, 0\leq t < + \infty
$$
and
$$
P\Bigg( X({t\wedge \tau^n}) = \xi + \int_{0}^{t} \tilde{b}(X(s))\mathbbm{1}_{\{s \leq \tau^n\}}\,ds
+ \varepsilon W({t \wedge \tau^n}); \,\,\,\forall\,\,\, 0\leq t < + \infty\Bigg)=1.
$$
\item $X$ has the strong Markov property, i.e. if we note $\tau:=\lim_{n
\rightarrow +\infty} \tau^n$ and $\tilde \tau$ is a stopping time of $X$ then,
conditional on $\tilde\tau<\tau$ and $X(\tilde\tau)=x$, the future $\{
X^+(t)=X(t+\tilde \tau) \colon t<\tau-\tilde\tau\}$ is independent of the past
$\{ X(s)\colon s\le \tilde\tau\}$ and identical in law to the process started at
$x$.
\end{itemize}
\end{defi}

We call $\tau$ the {\em explosion time} for $X$. Notice that the assumption of
continuity of $X$ in $\R^d\cup \{\infty\}$ implies that
$$
\tau = \inf \{ t > 0 : X(t) \notin \R^d \} \hspace{0.5cm}\mbox{ and }\hspace{0.5cm}X( \tau - )=X(\tau)=\infty\,\,\mbox{ on }\,\,\{ \tau < +\infty\}.
$$ We stipulate that $X(t)=\infty$ provided that $\tau \leq t < +\infty$ but we do not assume that $\lim_{t \to +\infty} X(t)$ exists when $\tau=+\infty$.

Notice that the assumption of finiteness of $\xi$ grants us $P(\tau > 0)=1$.
Also, if $P(\tau = +\infty)=1$ then we are left with the usual definition of
strong solution to the equation.

\begin{rem} It can be proved that if $\tilde{b} \in C^1(\R^d)$ then there
exists a unique solution of \eqref{eqexplosion0} up to an explosion time (see
\cite{KS, McKean}).
\end{rem}

\subsection{Freidlin-Wentzell estimates}

One of the most valuable tools in the study of perturbations by additive white
noise of an ODE is the Freidlin-Wentzell theory, whose main results we briefly
describe here.

Let $X^{x,\ve}$ be a solution to the SDE
$$
dX^{x,\ve} = \tilde{b}(X^{x,\ve})\,dt + \varepsilon\,dW
$$ with initial condition $x \in \R^d$, where $\tilde b$ is globally Lipschitz
with Lipschitz constant $K$.  
Fix $T > 0$ and let $P_x^{\ve,T}$ denote the law of $X^{x,\ve}$ on
$C([0,T],\R^d)$. Let us also consider $X^{x}$ the unique solution to the
deterministic equation
$$
\dot{X}(t)=\tilde b(X(t))
$$ with initial condition $x \in \R^d$.

\begin{teo}[Freidlin and Wentzell, \cite{FW}]
For each $x \in \R^d$ and $T>0$ the family $(P_x^{\ve,T})_{\ve > 0}$ satisfies a
large deviations principle on $C([0,T],\R^d)$ with scaling $\ve^{-2}$ and (good)
rate function $I^x_T$ given by
$$
I^x_T(\varphi)=\left \{\begin{array}{ll}\frac12 \int_0^T|\dot \varphi(s) - \tilde b(\varphi(s))|^2\, ds & \mbox{if $\varphi$ is absolutely continuous and $\varphi(0)=x$}\\
\\
+\infty & \mbox{otherwise}
                       \end{array} \right.
$$
\end{teo}

As a matter of fact, we need only the following weaker statement for our analysis: for every fixed $T > 0$ and $\delta >0$ there exist positive constants $C_1$ and $C_2$ depending on
$T$, $\delta$ and $K$ such that for all $0 < \varepsilon \leq 1$
\begin{equation}
 \label{grandes1}
\sup_{x \in \R^d} P\bigg( \sup_{t \in [0,T]} | X^{x,\,\varepsilon}(t) - X^{x}(t) | > \delta \bigg) \leq C_1 e^{-\frac{C_2}{\varepsilon^2}}.
\end{equation}

\subsection{Main results}

We now state the main results of the article. The first of them concerns the
explosion time of solutions to $\eqref{1.epsilon}$. In the following $P_u$
denotes the law of the solution to $\eqref{1.epsilon}$ up to the explosion time
$\tau^u_\ve$ with initial condition $u$. When the initial condition is clear
we often write $\tau_\ve$ instead of $\tau^u_\varepsilon$ to simplify the
notation.

\begin{teo}
\label{blowup1}
Let $U^{u,\ve}$ be a solution to \eqref{1.epsilon}. Then $P_u(\tau_\ve < \infty) = 1$.
\end{teo}

Let us notice that this result establishes a first difference in behavior with
respect to the deterministic system. While global solutions exist in the
deterministic equation, they do not for the stochastic one.

We then focus on establishing the order of magnitude and asymptotic
distribution of the explosion time for the different initial conditions $u \in
\R^d$. We deal first with initial conditions in the domain of explosion $D_ e$
and show the following result.

\begin{teo} \label{convergence} Given $\delta > 0$ and $u \in D_e$ we have
\begin{equation}\label{convergencia}
\lim_{\varepsilon \rightarrow 0} P_u ( |\tau_\varepsilon - \tau_0| > \delta ) = 0.
\end{equation} Moreover, the convergence is exponentially fast.
\end{teo}

This last theorem shows that for small $\varepsilon > 0$ the behavior of the
stochastic system  does not differ significantly from the deterministic one for
initial conditions in $D_e$. However, this is not the case for initial data in
the domain of attraction of the origin. Here is where important
differences appear and where characteristics associated with
metaestability are observed. In order to properly state the results achieved in
this matter, we need to introduce some notation.

For each $\varepsilon > 0$ we define
$$
\beta_{\varepsilon}= \inf \{ t \geq 0 : P_0 ( \tau_\ve > \beta_\varepsilon ) \leq e^{-1} \}
$$ which is well defined since $P_0 (\tau_\varepsilon < +\infty ) = 1$ for
every $\varepsilon > 0$.
We first show that the family $(\beta_\varepsilon)_{\varepsilon > 0}$ verifies
\begin{equation*}
\label{n4.2}
\lim_{\varepsilon \rightarrow 0} \varepsilon ^{2}\log\beta_{\varepsilon} = \Delta
\end{equation*}
with $\Delta := 2(\phi(\1) - \phi(0))$. In fact, we prove the stronger statement featured in the following theorem.

\begin{teo}\label{magnitude}
For each $u \in D_0$ and $\delta > 0$
\begin{equation*}\label{n4.3}
\lim_{\varepsilon \rightarrow +\infty} P_u \Big( e^{\frac{\Delta - \delta}{\varepsilon^2}} < \tau_\varepsilon < e^{\frac{\Delta + \delta}{\varepsilon^2}}\Big)=1,
\end{equation*} where the convergence is uniform over compact subsets of $D_0$.
\end{teo}
This theorem characterizes the asymptotic order of magnitude of the explosion
time for any initial condition $u \in D_0$. Regarding its distribution, we show
the asymptotic loss of memory in our last result.

\begin{teo} For each $u \in D_0$ and $t > 0$
\label{distribution}
\begin{equation*}\label{n4.4}
\lim_{\varepsilon \rightarrow 0} P_{u} (\tau_{\varepsilon} > t\beta_{\varepsilon})= e^{-t}
\end{equation*} where the convergence is uniform over compact subsets of $D_0$.
\end{teo}

\section{The deterministic system}
\label{deter}

\noindent Throughout this section we state some properties and study the behavior of solutions to \eqref{1.1}. This is carried out in \cite{AFBR} for solutions with nonnegative initial conditions. The purpose of this section is to extend the analysis in \cite{AFBR} to any arbitrary initial data $u\in \R^d$.

Let us start by noticing that equation \eqref{1.1} can be written as
\begin{equation*}\label{n3.1}
\dot{U}(t) = b(U(t))
\end{equation*} for $b=-\nabla \phi$ where $\phi$ is defined as
\begin{equation}
\label{lyapunov}
\phi(U) = \frac{1}{2} \langle AU , U \rangle - \frac{2}{h}\Big(\frac{\;\;\;|U_d^+|^{p+1}}{p+1} - \frac{\,\,{U_d}^2}{2}\Big).
\end{equation}
Here $A$ is as in \eqref{1.epsilon}-\eqref{A}. Notice that the potential $\phi$ has exactly two critical points:
$\1:=(1,\dots,1)$ and the origin. Both of them are hyperbolic. The origin is the only local minimum of $\phi$ while $\1$ is a saddle point. Our goal is to decompose $\R^d$ into distinct regions, each of them having different asymptotic characteristics under our system. To be able to accomplish such decomposition we need a few results concerning solutions to \eqref{1.1}. We begin with the following proposition.

\medskip
\begin{prop}\label{Lyapunov} Let $U=(U_1,\dots,U_d)$ be a solution to \eqref{1.1}. Then the application $t \mapsto \phi(U(t))$ is monotone decreasing.

\end{prop}

\begin{proof} Since $A$ is symmetric and $\dot{U} = -AU + \frac{2}{h}g(U_d)e_d$, a direct calculation shows that
$$
\frac{d\phi\big(U(t)\big)}{dt} = \langle \dot{U}(t) ,\, AU(t) \rangle - \frac{2}{h}g\big(U_d(t)\big)\dot{U}_d(t)= - |\dot{U}(t)|^2 \leq 0.
$$
\qedhere
\end{proof}
Next we show that solutions to \eqref{1.1} satisfy a Maximum Principle.

\begin{lema}[\textbf{Maximum Principle}]\label{prinmax} Let $U=(U_1,\dots,U_d)$ be a solution to
\eqref{1.1}.
Then $U$ satisfies
\begin{equation}\label{pmaximo} \max_{k=1,\dots,d}
|U_k (t)| \leq \max \{ \max_{k=1,\dots,d} |U_k(0)| , \max_{0\leq s \leq t} U_d(s)\}
\end{equation}
\end{lema}

\begin{proof} We prove first that
\begin{equation}\label{pmaximo1}
\max_{k=1,\dots,d}
|U_k (t)| \leq \max \{ \max_{k=1,\dots,d} |U_k(0)| , \max_{0\leq s \leq t} |U_d(s)|\}
\end{equation} and then we check that if \eqref{pmaximo1} holds then
$$
\max \{ \max_{k=1,\dots,d} |U_k(0)| , \max_{0\leq s \leq t} |U_d(s)|\}=\max \{ \max_{k=1,\dots,d} |U_k(0)| , \max_{0\leq s \leq t} U_d(s)\}
$$
which allows us to conclude \eqref{pmaximo}. Let $j$ be the node that maximizes
$\max_{0 \leq s \leq t} |U_j(s)|$. Let us observe that if $j=d$ then
\eqref{pmaximo1} is immediately verified. Hence, we can assume that $1\leq j <
d$. Consider $t_0 = \min\{ t' \in [\,0,t\,] : \max_{0 \leq s \leq t}
|U_j(s)| = |U_j(t')|\}$, the first time in which the maximum is attained. Note that $|U_j(t_0)|=\max_{k=1,\dots, d}( \max_{0 \leq s \leq t} |U_k(s)|)$. If $t_0=0$ then
$$
\max_{k=1,\dots,d} |U_k(0)| \geq |U_j(t_0)| = \max_{k=1,\dots,d} \Big(
\max_{0 \leq s \leq t} |U_k(s)| \Big) \geq \max_{k=1,\dots,d} |U_k (t)|
$$ and we get \eqref{pmaximo1}. If $t_0 > 0$ we must consider two cases: $U_j(t_0) \geq 0$ and
$U_j(t_0) < 0$. If $U_j(t_0) \geq 0$ then by definition of $t_0$ we get that
$U_j(t_0) \geq U_j(s)$ for all $0 \leq s \leq t$. From this it follows that
$U'_j(t_0)\geq 0$. On the other hand, the choice of $j$ guarantees that
$U_j(t_0) \geq U_k(t_0)$ for all $k=1,\dots,d$. This implies that
$$
U'_j(t_0) = \frac{1}{h^2} \big( (U_{j+1}(t_0) - U_j(t_0)) + (U_{j-1}(t_0) - U_j(t_0) )\big) \leq 0 \qquad \text{if}\,\, 1 < j < d
$$ and
$$
U'_1(t_0) = \frac{2}{h^2} ( U_2(t_0) - U_1(t_0) ) \leq 0 \qquad \text{if}\,\, j=1.
$$
In any case we conclude that $U'_j(t_0)=0$ and, in particular, that
$U_{j+1}(t_0)= U_j(t_0)$. We conclude that $|U_{j+1}(t_0)|=\max_{k=1,\dots, d} (
\max_{0 \leq s \leq t}
|U_k(s)|)$ which allows us to repeat the same argument, now for $j+1$ instead
of $j$. Thus, an inductive procedure eventually yields that
$U_{d}(t_0)=U_j(t_0)$. From here we obtain \eqref{pmaximo1} if $U_j(t_0) \geq
0$. The case $U_j(t_0) < 0$ is analogous. To conclude \eqref{pmaximo} we notice
that if $t_1 = \min\{ t' \in [\,0,t\,] : \max_{0 \leq s \leq t}
|U_d(s)| = |U_d(t')|\} > 0$ then $U_d(t_1)\geq 0$ because, otherwise, from
\eqref{1.1} and \eqref{pmaximo1} we get that $U'_d(t_1)>0$ which contradicts the
definition of $t_1$.
\end{proof}

\medskip
As a consequence of the Maximum Principle we have the following
characterization of globally defined solutions to \eqref{1.1}.

\medskip
\begin{lema}\label{lemaacot} Let $U$ be a globally defined solution to
\eqref{1.1}. Then $U$ is bounded.
\end{lema}

\begin{proof} Let us suppose that $U$ is not bounded. Then by the Maximum
Principle we obtain that $\max_{0\leq s\leq t} |U_d(s)| \rightarrow +\infty$ as
$t \rightarrow +\infty$.

\noindent \textbf{1}. Given $M > 0$ we define $t_M := \inf \{ t \geq
0 : |U_d(t)| > M \}$. From this definition it follows that $|U_d(t_M)| \geq M$
and that $|U_d(t_M)|=\max_{0 \leq s \leq t_M }
|U_d(s)|$. If $M > \max_{k=1,\dots,d} |U_k(0)|$ then $t_M > 0$ and by the
Maximum Principle we have $|U_{d-1}(t_M)| \leq  U_d (t_M)$. This gives us the
inequality
\begin{equation*}\label{ndesig}
U'_d(t_M) \geq \frac{2}{h} U_d^p(t_M) - \Big( \frac{4}{h^2} + \frac{2}{h}\Big) U_d(t_M).
\end{equation*}

\bigskip
\noindent \textbf{2}. From here it is easy to see that if $M$ is large enough
we have that $U_d : [t_M , +\infty) \longrightarrow \R$ is monotone increasing.
This implies that for $t \geq t_M$ we have $U_d(t)=\max_{0\leq s \leq t}
|U_d(s)|\geq M $ and, as a consequence, that $U'_d(t) \geq \frac{2}{h} U_d^p(t)
- \big(\frac{4}{h^2} + \frac{2}{h}\big) U_d(t)$. If $M$ is taken large enough
then $U$ verifies $U'_d(t) \geq \frac{1}{h} U_d^p(t)$ for $t \geq t_M$ and,
therefore, cannot be globally defined. This is a contradiction which implies
that $U$ must be bounded.
\end{proof}

From the previous lemma and the fact that \eqref{1.1} admits the Lyapunov
functional \eqref{lyapunov} we obtain the following corollary.

\begin{cor}\label{corcomportamiento} Let $U$ be a solution to \eqref{1.1}. Then
either $U$ explodes in finite time or is globally defined and converges to a
stationary solution as $t \rightarrow +\infty$.
\end{cor}

With this result at our disposal we can obtain the following theorem, whose
proof is in \cite{AFBR}.

\begin{teo} {\ }
\label{divisionposta}
\begin{enumerate}
\item Equation \eqref{1.1} has exactly two equilibriums $U\equiv 0$ and $U\equiv
1$. The first one is stable and the second one is unstable.

\item Let $u$  be a nonnegative initial datum such that $U^u$ is globally
defined and $\lim_{t \rightarrow +\infty} U^u(t)= \1$. Then
\begin{itemize}
\item $0\le v \lneq u \Longrightarrow U^v$ is globally defined and $
\displaystyle \lim_{t \rightarrow +\infty} U^v(t)= 0.$
\item $u\lneq v \Longrightarrow U^v$ explodes in finite time.
\end{itemize}
\item Consider $\lambda > 0$ and a nonnegative initial condition $u$. Then
there exists $\lambda_c > 0$ such that
\begin{enumerate}
\item $\lambda < \lambda_c \Longrightarrow U^{\lambda u}$ is globally defined
and $\lim_{t \rightarrow +\infty} U^{\lambda u} (t)=0$
\item $\lambda_c < \lambda \Longrightarrow U^{\lambda u}
    \,\,\,\text{explodes in finite time}$
\item $ \lambda = \lambda_c \Longrightarrow U^{\lambda u}$ is globally defined
and $\lim_{t \rightarrow +\infty}U^{\lambda u}(t)= \1$.
\end{enumerate}
\end{enumerate}
\end{teo}

\begin{figure}
	\centering
	\includegraphics[width=8cm]{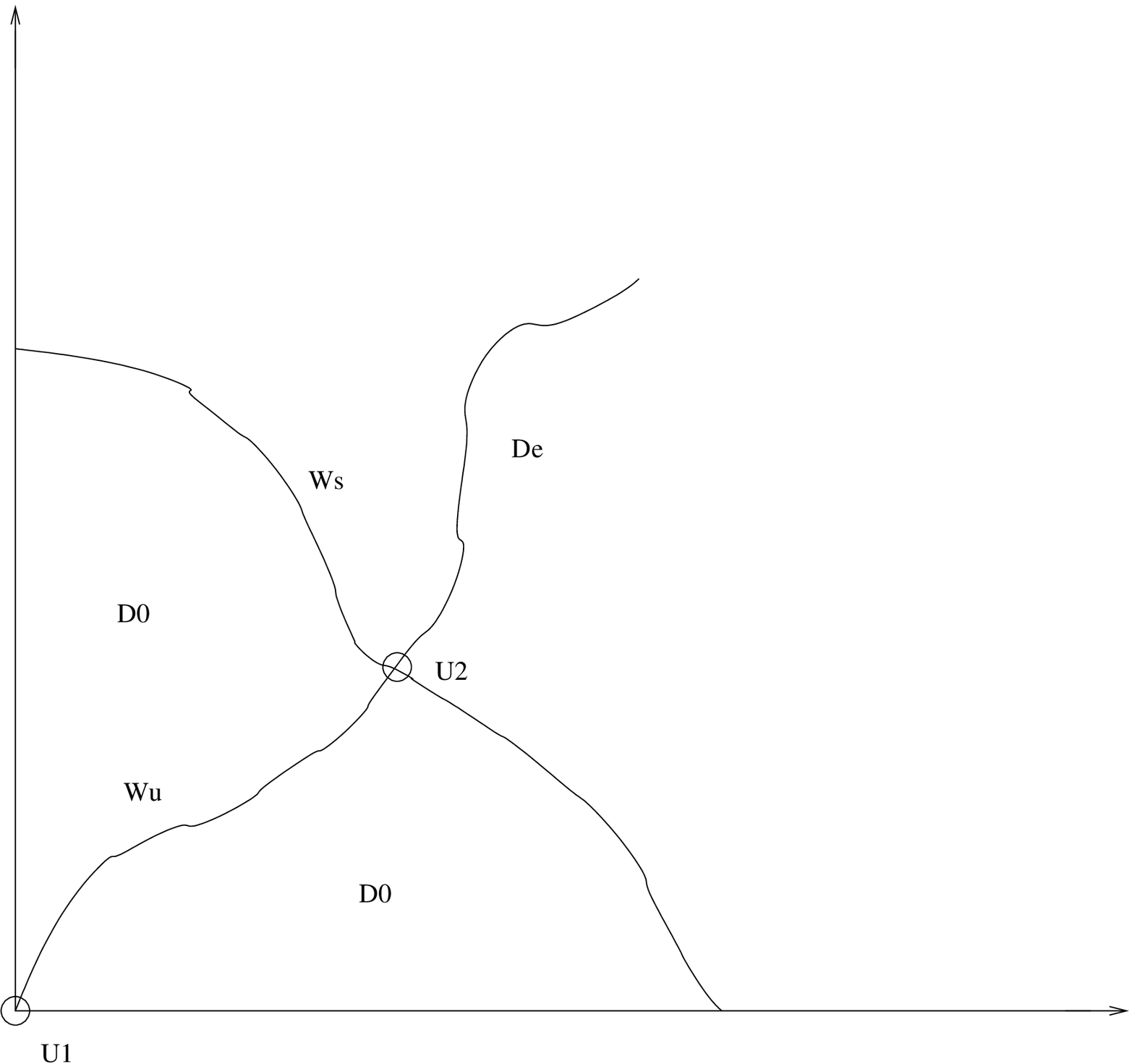}
	\caption{The phase diagram of equation \eqref{1.1}.}
	\label{fig:inclination}
\end{figure}

This results allow us to give a good description of the behavior of the
deterministic system $U$ for the different initial conditions $u \in \R^d$.
Indeed, we have a decomposition
\begin{equation*}
\R^d = D_0 \cup \mathcal{W}^s_1 \cup D_e
\end{equation*} where $D_0$ denotes the stable manifold of the origin,
$\mathcal{W}^s_{1}$ is the stable manifold of $\1:=(1,\dots,1)$ and $D_e$ is the
domain of explosion, i.e., if $u \in D_e$ then $U^{u}$ explodes in finite time.
The sets $D_0$ and $D_e$ are open in $\R^d$. The origin is an asymptotically
stable equilibrium of the system. $\mathcal{W}^s_{1}$ is a manifold
of codimension one. Also $\1$ admits an unstable manifold
of dimension one which we shall note by $\mathcal W^u_\1$. This unstable
manifold is contained in $\R^d_+$, has nonempty intersection with both $D_0$ and
$D_e$ and joins $\1$ with the origin.
An illustration of this decomposition is given in Figure \ref{fig:inclination}
for the $2$-dimensional case.

\psfrag{U1}{\vspace{-155pt}$U\equiv 0$}
\psfrag{U2}{\vspace{-45pt}$U \equiv 1$}
\psfrag{De}{\hspace{60pt}\vspace{30pt}$D_e$}
\psfrag{D0}{$D_0$}
\psfrag{Wu}{$\mathcal W_1^u$}
\psfrag{Ws}{$\mathcal W_1^s$}

\section{Explosions in the stochastic model}
\label{estoc.exp}

In this section we focus on proving that solutions to \eqref{1.epsilon} blow up
in finite time with probability one for any initial condition $u \in \R^d$ and
every $\ve >0$. The idea is to show that, conditioned on non-explosion, the
system is guaranteed to enter a specific region of space in which we can prove
that explosion occurs with total probability. From this we can conclude that
non-explosion must happen with zero probability. We do this by comparison with
an adequate Ornstein-Uhlenbeck process.

\begin{proof}[Proof of Theorem \ref{blowup1}]
Let $Y^{y,\,\varepsilon}$ be the solution to
\begin{equation*}
\label{OU}
dY^{y,\,\varepsilon}=-\Big(AY^{y,\,\varepsilon} + \frac{2}{h} Y^{y,\,\varepsilon}_d e_d\Big)\,dt + \varepsilon dW
\end{equation*}
with initial condition $Y^{y,\ve}(0)=y$. Notice that the drift term is linear,
and  given by a negative definite matrix. Hence, $Y^{y,\ve}$ is in fact a
$d$-dimensional Ornstein-Uhlenbeck process which admits an invariant
distribution supported in $\R^d$. We also have convergence to this equilibrium
measure for any initial distribution and therefore the hitting time of
$Y^{y,\ve}$ of any open set is finite almost surely.

On the other hand, since the drift term of \eqref{OU} is smaller or equal than
$b$ we can apply the stochastic comparison principle to obtain that
$U^{u,\ve}(t) \ge Y^{y,\ve}$ holds a.s. as long as $U^{u,\ve}$ is finite, if
$u\ge y$. From here, the result follows applying the following lemma and the
strong Markov property.
\end{proof}

\begin{lema}
Consider the set
$$
\Theta^M:=\{ y \in \R^d : y_k \geq 0 \,\,\text{for all}\,\,0\leq k \leq d-1\,,\, y_d \geq M \},
$$
then we have
$$
\lim_{M\to \infty }\inf_{y\in \Theta^M} P_y(\tau_\ve <\infty) =1.
$$
\end{lema}

\begin{proof}
Consider the auxiliary process $Z^{y,\,\varepsilon}:= U^{y,\,\varepsilon} - \ve
W$. Notice that this process verifies the random differential equation
\begin{equation*}
 dZ^{y,\ve}=b(Z^{y,\ve} +\ve W)dt, \quad Z^{y,\ve}(0)=y.
\end{equation*}
Let us also observe that $Z^{y,\ve}$ has the same explosion time as
$U^{y,\ve}$. For each $k \in \N$ let us define the set $A_k:=\{\sup_{0\leq t
\leq 1} |W_d(t)|\le k\}$. On $A_k$ we have that $Z^{y,\ve}$ verifies the
inequality
\begin{equation}
\label{ineq.z}
\frac{dZ^{y,\,\varepsilon}}{dt} \ge -AZ^{y,\,\varepsilon} - \frac{4}{h^2}\ve k \sum e_i+ \frac{2}{h} ((Z^{y,\,\varepsilon}_d - \varepsilon k)_+^p - Z^{y,\,\varepsilon}_d - \varepsilon k) e_d.
\end{equation}
Observe that \eqref{ineq.z} can be written as
$$
\frac{dZ^{y,\,\varepsilon}}{dt} \ge QZ^{y,\,\varepsilon} + q + (Z^{y,\,\varepsilon}_d - \varepsilon k)_+^p \ge QZ^{y,\,\varepsilon} + q ,
$$
where $Q\in\R^{d\times d}$ and $q\in\R^d$ both depend on $\ve, h$ and $k$, but
not on $M$. This allows us to conclude the inequality $|Z^{y,\,\varepsilon}| \le
(M + |q|){\rm exp}(|Q|)$ for all $0\le t\le 1$. In particular, for the last
coordinate we get
$$
\left\{\begin{array}{ll}
\frac{dZ_d^{y,\,\varepsilon}}{dt} \ge -\alpha_1 M + \alpha_2
(Z^{y,\,\varepsilon}_d)^p & \mbox{ if }\,\,0\le t \le 1\\
\\
Z_d^{y,\,\varepsilon}(0)\ge M
\end{array}
\right.
$$
for constants $\alpha_1, \alpha_2$ which do not depend on $M$. It is a
straightforward calculation to check that solutions to this one-dimensional
inequality blow-up in a finite time that converges to zero as \mbox{$M\to
+\infty$.} Therefore, for each $k \in \N$ there exists $M_k$ such that
$P(A_k)\leq \inf_{y \in \Theta^{M}}  P_y(\tau_\ve <\infty)$ for all $M \geq
M_k$. Since $\lim_{k \rightarrow +\infty} P(A_k)=1$, this concludes the proof.
\end{proof}

\section{Convergence of $\tau^u_\varepsilon$ for initial conditions in $D_e$}
\label{estoc.conv}

This section is devoted to prove that for initial data in the domain of
explosion of the deterministic system, the explosion time is of order one and,
moreover, as $\ve \to 0$ converges to the explosion time of the deterministic
system. Observe that do to the lack of boundedness this result do not follow
from standard perturbation arguments for dynamical system (deterministic or
stochastic). We first introduce the truncations of the drift that we use here to
prove one of the bounds and we are going to make more profit of them in Section
\ref{meta} when we deal with initial data in the domain of attraction of the
origin.

\subsection{Truncations of the potential and localization}\label{trunca}

\noindent The large deviations principle originally formulated by Freidlin and
Wentzell for solutions of stochastic differential equations like
\eqref{eqexplosion0} require a global Lipschitz condition on the drift term
$\tilde b$. While this condition is met on the classic double-well potential
model, it is not in our case. As a consequence, we cannot apply such estimates
to our system directly. Nonetheless, the use of localization techniques helps us
to solve this problem and allows us to take advantage of the theory developed by
Freidlin and Wentzell despite the fact that our drift term is not globally
Lipschitz. In the following lines we give details about the localization
procedure to be employed in the study of our system.

\medskip

For every $n \in \N$ let $G_n : \R \longrightarrow \R$ be of class $C^2$ such that
\begin{equation*}\label{gtruncada}
G_n(u) = \left\{\begin{array}{ll}
\frac{|u^+|^{p+1}}{p+1} - \frac{u^2}{2} &\,\,\text{if}\,\,|u| \leq n\\
0 &\,\,\text{if}\,\,|u| \geq 2n.
\end{array}\right.
\end{equation*}
We consider then the family $\big(\phi^{n}\big)_{n \in \N}$ of potentials over $\R^d$ given by
\begin{equation*}
\phi^{n}(u)= \frac{1}{2}\langle Au, u \rangle - \frac{2}{h} G_n(u_d).
\end{equation*} This family satisfies the following properties:

\begin{itemize}
\item [(i)]For every $n \in \N$ the potential $\phi^{n}$ is of class $C^2$ and $b^{n}=- \nabla \phi^{n}$ is globally Lipschitz.
\item [(ii)]For $n \leq m \in \N$ we have $b^{n}\equiv b^{m}$ over the region $\Pi^n = \{ u \in \R^d : |u_d| < n \}.$
\item [(iii)]For every $n \in \N$ we have $\liminf_{|u| \rightarrow +\infty} \frac{\phi^{n}}{|u|} > 0.$
\end{itemize}
Since $b^{n}$ is globally Lipschitz, for each $u \in \R^d$ there exists a unique solution to the ordinary differential equation
\begin{equation*}
\dot{U}^{n,\,u} = b^{n}(U^{n,\,u})
\end{equation*} with initial condition $u$. Such solution is globally defined
and describes the same trajectory as the solution to \eqref{1.1} starting at $u$
until the escape from $\Pi^n$.
In the same way, for each $x \in \R^d$ and $\varepsilon > 0$ there exists a
unique global solution to the stochastic differential equation
\begin{equation}\label{eqtruncada}
dU^{n,\,u,\,\varepsilon} = b^{n}\big(U^{n,\,u,\,\varepsilon}\big)\,dt + \varepsilon\,dW
\end{equation} with initial condition $u$.

As before we use $U^{n,\, u}$ for $U^{n,\, u,\, 0}$. Since $b^{n}$ coincides
with $b$ over the ball $B_n(0)$ of radius $n$ centered at the origin, if we
write
\begin{equation*}
\label{tiempostau}
\tau^{n,\, u}_\varepsilon = \inf \{ t \geq 0 : |U^{n,\,u,\,\varepsilon}(t)|
\geq n \}, \qquad \tau_\varepsilon^u:=\lim_{n \rightarrow +\infty} \tau^{n,\,
u}_\varepsilon,
\end{equation*}
then for $t < \tau_\varepsilon^u$ we have that $U^{u,\,\varepsilon}(t) := \lim_{n \rightarrow +\infty} U^{n,\,u,\,\varepsilon}(t)$ is a solution to
\begin{equation}\label{eq}
dU^{u,\,\varepsilon} = b(U^{u,\,\varepsilon})\,dt + \varepsilon\,dW
\end{equation} until the explosion time $\tau^u_\varepsilon$ with initial condition $u$. Moreover, if we define the stopping times
\begin{equation*}\label{tiempospi}
\pi^{n,\,u}_\varepsilon = \inf \{ t \geq 0 : U^{n,\,u,\,\varepsilon}(t) \notin \Pi^n \},
\end{equation*} it can be seen that (ii) implies that
$$
\tau^u_\varepsilon= \lim_{n \rightarrow +\infty} \pi^{n,\,u}_\varepsilon
$$ and that $U^{u,\,\varepsilon}$ coincides with the process
$U^{n,\,u,\,\varepsilon}$ until the escape from $\Pi^n$. On the other hand,
(i) guarantees that for each $n \in \N$ and $u \in \R^d$ the family
$\big(U^{n,\,u,\,\varepsilon}\big)_{\varepsilon > 0}$ satisfies a large
deviations principle. Finally, from (iii) we get that there is an unique
invariant probability measure for the process $U^{n,\,\varepsilon}$ for each
$\varepsilon > 0$ given by the formula
\begin{equation*}\label{ninvariante}
\displaystyle \mu_{\varepsilon}^{n}(A):= \frac{1}{Z_\varepsilon^{n}} \int_{A} e^{-\frac{2}{\ve^2}\phi^{n}(u)}\,du, \qquad A \in \mathcal{B}(\R^d)
\end{equation*}where $Z_\varepsilon^{n}= \int_{\R^d}
e^{-\frac{2\phi^{n}(u)}{\varepsilon^2}}\,du.$
Hereafter, when we refer to the solution of \eqref{eq} we mean the solution
constructed in this particular way.

\subsection{Proof of Theorem \ref{convergence}}

We split the proof of Theorem \ref{convergence} in two parts, the first one is
immediate from the continuity of the solutions of \eqref{1.epsilon} with respect
to $\ve$ in intervals where the deterministic solution is bounded.

\begin{prop} For any fixed $\delta > 0$ and $u \in D_e$ we have
\begin{equation*}\label{convergenciainferior}
\lim_{\varepsilon \rightarrow 0} P_u ( \tau_\varepsilon < \tau_0 - \delta ) = 0.
\end{equation*}
\end{prop}

\begin{proof} We may assume that $\tau^u_0 > \delta$ since the proof is trivial
otherwise. Now, as the deterministic system $U^{u}$ is defined up until
$\tau^u_0$, if we take $M:= \sup_{0 \leq t \le \tau^u_0 - \delta} |U^{u}_t| <
+\infty$ then $\tau^u_\varepsilon < \tau_u - \delta$ implies that
$$
\sup_{0 \le t \le \tau_u-\delta} \big| U^{2M,\,u,\,\varepsilon}(t)- U^{2M,\,u}(t) \big| > 1.
$$ By  \eqref{grandes1} we get \eqref{convergenciainferior}.
\end{proof}

\begin{prop}\label{convsup} For any $\delta > 0$ and $u \in D_e$ we have
\begin{equation*}\label{convergenciasuperior}
\lim_{\varepsilon \rightarrow 0} P_u ( \tau_\varepsilon > \tau_0 + \delta ) = 0.
\end{equation*}
Moreover, the convergence is uniform over compact subsets of $D_e$.
\end{prop}

\begin{proof}Fix $\delta > 0$, $\mathcal{K}$ a compact set contained in $D_e$
and let $Y^{u}$ be the solution to the ordinary differential equation
$$
\dot{Y}^{u}=-\Big(AY^{u} + \frac{2}{h} Y^{u,\,\varepsilon}_d e_d\Big)
$$ with initial condition $u \in \mathcal{K}$. By the Comparison principle we
have that $U^{u} \geq Y^{u}$ for as long as $U^{u}$ is defined. Since $Y^{u}$ is
the solution to a linear system of ordinary differential equations whose
associated matrix is symmetric and negative definite, we get that there exists
$\rho_{\mathcal{K}} \in \R$ such that for all $u \in \mathcal{K}$ every
coordinate of $U^{u}$ remains bounded from below by $\rho_{\mathcal{K}} + 1$ up
until $\tau^u_0$.
If for $\rho \in \R$ and $M > 0$ we write
$$
\Theta_{\rho}^M:=\{ y \in \R^d : y_k \geq \rho \,\,\text{for all}\,\,0\leq k \leq d-1\,,\, y_d \geq M \}
$$ then by the Maximum Principle and the previous statement we have that $T_u:=
\inf \{ t \geq 0 : U^{u}_t \in \Theta_{\rho_{\mathcal{K}}+1}^{M+1}\}$ is finite.
Moreover, as $U^{M+2,\,u}$ agrees with $U^{u}$ until the escape from
$\Pi_{M+2}$, we obtain the expression $T_u= \inf \{ t \geq 0 : U^{M+2,\,u}_t \in
\Theta_{\rho_{\mathcal{K}+1}}^{M+1}\}$. Taking $T_{\mathcal{K}}:=\sup_{u \in
\mathcal{K}} T_u <+\infty$ we may compute
\begin{align*}
P_u\big(\tau_{\varepsilon}(\Theta_{\rho_{\mathcal{K}}}^M) > T_u\big) &\leq
P_u\big( \pi^{M+2}_\varepsilon \wedge
\tau_{\varepsilon}(\Theta_{\rho_{\mathcal{K}}}^M) > T_u \big) +
P_u\big(\pi^{M+2}_\varepsilon \leq
T_u\,,\,\tau_{\varepsilon}(\Theta_{\rho_{\mathcal{K}}}^M) > T_u \big)\\
\\
& \leq 2 P_u \Big(\sup_{0\le t \le T_u} |U^{M+2,\,\varepsilon}(t) - U^{M+2}(t)| > 1 \Big)\\
\\
& \leq 2 P_u \Big(\sup_{0\le t \le T_{\mathcal{K}}} |U^{M+2,\,\varepsilon}(t) - U^{M+2}(t)| > 1 \Big),
\end{align*}
from which by \eqref{grandes1} we obtain 
\begin{equation}\label{convergenciaeq1}
\lim_{\varepsilon \rightarrow 0} \sup_{u \in \mathcal{K}} P_u\big(\tau_{\varepsilon}(\Theta_{\rho_{\mathcal{K}}}^M) > T_u\big) = 0.
\end{equation}
On the other hand, by the strong Markov property for $U^{u,\varepsilon}$ we get
$$
P_u \big( \tau_\varepsilon > \tau_0 + \delta \big) \leq P_u \big(
\tau_\varepsilon > T_u + \delta \big) \leq \sup_{y \in
\Theta_{\rho_{\mathcal{K}}}^M} P_y ( \tau_\varepsilon > \delta) + \sup_{u \in
\mathcal{K}} P_u\big(\tau_{\varepsilon}(\Theta_{\rho_{\mathcal{K}}}^M) >
T_u\big).
$$
Taking into consideration \eqref{convergenciaeq1}, in order to finish the
proof we only need to show that the first term on the right hand side tends to
zero as $\varepsilon \rightarrow 0$ for an adequate choice of $M$. To see this
we consider for each $\varepsilon > 0$ and $y \in \Theta_{\rho_{\mathcal{K}}}^M$
the processes $Y^{y,\,\varepsilon}$ and $Z^{y,\,\varepsilon}$ defined by
$$
dY^{y,\,\varepsilon}=-\Big(AY^{y,\,\varepsilon} + \frac{2}{h}
Y^{y,\,\varepsilon}_d e_d\Big)\,dt + \varepsilon dW,
$$
and $Z^{y,\,\varepsilon}:= U^{y,\,\varepsilon} - Y^{y,\,\varepsilon}$,
respectively. Notice that since $Y^{y,\,\varepsilon}$ is globally defined and
both $U^{y,\,\varepsilon}$ and $Z^{y,\,\varepsilon}$ have the same explosion
time. Also note that $Z^{y,\,\varepsilon}$ satisfies the random differential
equation
$$
dZ^{y,\,\varepsilon}=-\Big(AZ^{y,\,\varepsilon} + \frac{2}{h}\Big(
\Big[\big(U_d^{y,\ve}\big)^{+}\Big]^p - Z^{y,\,\varepsilon}_d\Big)e_d\Big)\,dt.
$$
The continuity of trajectories allows us to use the Fundamental Theorem of
Calculus to show that almost surely $Z^{y,\,\varepsilon}(\omega)$ is a solution
to the ordinary differential equation
\begin{equation}\label{rde1}
\dot{Z}^{y,\,\varepsilon}(t)(\omega) = -AZ^{y,\,\varepsilon}(\omega) +
\frac{2}{h}\Big( \Big[\big(U_d^{y,\ve}\big)^{+}\Big]^p(\omega) -
Z^{y,\,\varepsilon}_d(\omega)\Big)e_d.
\end{equation}
For each $y \in \Theta_{\rho_{\mathcal{K}}}^M$ and $\varepsilon > 0$ let
$\Omega^y_\varepsilon$ be a set of probability one in which \eqref{rde1} holds.
Notice that for every $\omega \in \Omega^y_\varepsilon$ we have the inequality
\begin{equation*}\label{rde2}
\dot{Z}^{y,\,\varepsilon}(\omega) \geq -AZ^{y,\,\varepsilon}(\omega) - \frac{2}{h}Z^{y,\,\varepsilon}_d(\omega)e_d.
\end{equation*}
Using the Comparison Principle we conclude that
$Z^{y,\,\varepsilon}(\omega)\geq 0 $ for every $\omega \in \Omega^y_\varepsilon$
and, therefore, that the inequality $U^{y,\,\varepsilon}(\omega) \geq
Y^{y,\,\varepsilon}(\omega)$ holds for as long as $U^{y,\,\varepsilon}(\omega)$
is defined.

For each $y \in \Theta_{\rho_{\mathcal{K}}}^M$ and $\varepsilon > 0$ let us
also consider the set
$$
\tilde{\Omega}^y_\varepsilon = \Big\{ \omega \in \Omega : \sup_{0\le t\le
\delta} | Y^{y,\,\varepsilon}(\omega,t)-  Y^{y}(\omega,t)| \leq 1 \,,\, \sup_{0
\leq t \leq \delta} |\varepsilon W(\omega,t)| \leq 1 \Big\}.
$$ Note that  $\lim_{\varepsilon \rightarrow 0} \inf_{y \in \Theta_{\rho_{\mathcal{K}}}^M} P( \tilde{\Omega}^y_\varepsilon ) = 1$. Our goal is to show that if $M$ is chosen adequately then for fixed $y \in \Theta_{\rho_{\mathcal{K}}}^M$ the trajectory $U^{y,\,\varepsilon}(\omega)$ explodes before time $\delta$ for all $\omega \in \Omega^y_\varepsilon \cap \tilde{\Omega}^y_\varepsilon$.
From this we get that
$$
\inf_{y \in \Theta_{\rho_{\mathcal{K}}}^M} P(\tilde{\Omega}^y_\varepsilon) =
\inf_{y \in \Theta_{\rho_{\mathcal{K}}}^M} P( \Omega^y_\varepsilon \cap
\tilde{\Omega}^y_\varepsilon) \leq \inf_{y \in \Theta_{\rho_{\mathcal{K}}}^M}
P_y (\tau_\varepsilon \leq \delta ).
$$ and by letting $\varepsilon \rightarrow 0$ we conclude the result.

So let us take $y \in  \Theta_{\rho_{\mathcal{K}}}^M$, $\omega \in
\Omega^y_\varepsilon \cap \tilde{\Omega}_\varepsilon$ and suppose that
$U^{y,\,\varepsilon}(\omega)$ is defined in the interval $[0,\delta]$.
Notice that since $\omega \in \Omega^y_\varepsilon \cap
\tilde{\Omega}_\varepsilon$ then the $(d-1)$-th coordinate of
$Y^{y,\,\varepsilon}(\omega,t)$ is bounded from below by $\rho_{\mathcal{K}} -
1$ for $t \in [0,\delta]$. By comparison we know that the $(d-1)$-th coordinate
of $U^{y,\,\varepsilon}_t(\omega,t)$ is bounded from below by
$\rho_{\mathcal{K}} - 1$ as well.

From here we deduce that the last coordinate of $U^{y,\varepsilon}(\omega)$
verifies the integral equation
\begin{equation*}
{U}^{y,\,\varepsilon}_d(\omega,t) \geq {U}^{y,\,\varepsilon}_d(\omega,s) + \int_s^t \frac{2}{h²}\,\Big(- U^{y,\,\varepsilon}_d(\omega,r) +\rho_{\mathcal{K}} -1  + hg\big( U^{y,\,\varepsilon}_d(\omega,r)\big) \Big)\,dr - 1
\end{equation*}for $s < t \in [0,\delta]$. We can take $M \in \N$ large enough to guarantee that there exists a constant $\alpha > 0$ such that for all $m\geq M$ we have
\begin{equation*}\label{lcotafea}
\frac{2}{h²}\,\big(-m +\rho_{\mathcal{K}} -1  + hg(m) \big) \geq \alpha m^p.
\end{equation*} If we recall that $U^{y,\,\varepsilon}_d(\omega,0) \geq M$ then our selection of $M$ implies that
\begin{equation*}
{U}^{y,\,\varepsilon}_d(\omega,t) \geq M-1 + \alpha\int_0^t \big( U^{y,\,\varepsilon}_d(\omega,u)\big)^p\,du
\end{equation*}
for all $t \in [0,\delta]$. But if this inequality holds and $M$ is large
enough, one can check that $U^{y,\,\varepsilon}(\omega)$ explodes before time
$\delta$, which contradicts our assumptions. Therefore, if $y \in 
\Theta_{\rho_{\mathcal{K}}}^M$ and $\omega \in \Omega^y_\varepsilon \cap
\tilde{\Omega}_\varepsilon$ then $U^{y,\,\varepsilon}(\omega)$ explodes before
time $\delta$ and this fact concludes our proof.
\end{proof}

Combining these two propositions we get Theorem \ref{convergence}. Observe that
the bounds obtained decay to zero exponentially fast due to Proposition
\eqref{grandes1}. 
\section{Metastable behavior for initial conditions in $D_0$}
\label{meta}

Finally we focus on initial data in $D_0$, where the metastability phenomenon
can be appreciated. We start with the construction of an auxiliary domain that
contains the origin and such that the exit time from this domain is
asymptotically equivalent to the explosion time.

\subsection{Construction of an auxiliary domain}\label{dominioauxiliar}

\noindent In order to proceed with our analysis of the explosion time we must
first construct an auxiliary bounded domain. The purpose behind this
construction is to reduce our problem to a simpler one, the escape from this
domain. This is easier because we may assume that the drift coefficient $b$ is
globally Lipschitz, as the escape only depends on the behavior of the system
while it remains inside a bounded region. In this case, large deviations
estimates as the ones proved by Freidlin and Wentzell apply. We need a bounded
domain $G$ which verifies the following properties:
\begin{enumerate}
\item [$(1)$] $G$ is bounded, contains $\1$ and the origin.
\item [$(2)$] There exists $c > 0$ such that $B_c(0) \subseteq G$ and for all
$y \in B_c(0)$ the system $U^{y}$ is globally defined and tends to zero without
escaping $G$.
\item [$(3)$] The border of $G$ can be decomposed in two parts: $\partial^1$
and $\partial G \setminus \partial^1 $. The region of the border $\partial^1$ is
closed and satisfies
$\min_{u \in \partial G} \phi(u) = \min_{u \in \partial^1} \phi(u)$ and
$$\inf_{u \in \partial G \setminus \partial^1 } \phi(u) > \min_{u \in \partial
G} \phi(u).$$
\item [$(4)$] For all $y \in \partial^1$ the deterministic system $U^{y}$ explodes in finite time.

\end{enumerate}

The domain $G$ can be constructed as follows. Let us consider the value of
$\phi$ at the saddle point $\1$, $\phi(\1)=-1/(p+1) + 1/2>0=\phi(0)$ and $c>0$
such that $\phi(u)<\phi(\1)$ for $u\in B_c(0)$.

\begin{figure}
\begin{center}
 \includegraphics[width=7cm]{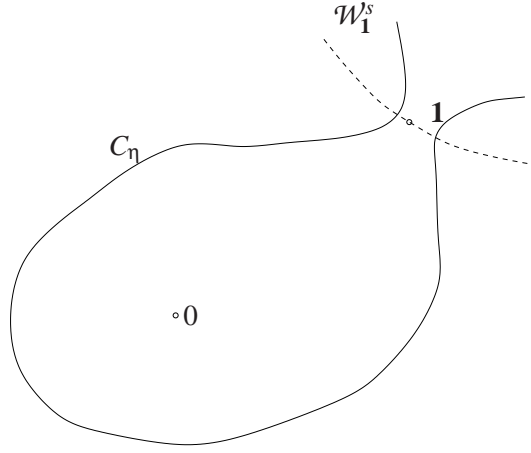}
\end{center}
\label{Gaux}
\caption{The level curve $C_\eta$ and the stable manifold of $\1$.}
\end{figure}

For each point $u\in \partial B_c(0)$ consider the ray $R_u:=\{ \lambda u :
\lambda > 0\}$. Since the vector $\1$ is not tangent to $\W_\1^s$ at $\1$, we
may take a sufficiently small neighborhood $V$ of $c\1$ such that for all $u \in
V\cap \partial B_c(0)$ the ray $R_u$ intersects $\W_\1^s \cap (\R_{> 0})^d$. For
such $V$ we may then define $\bar{\lambda}_u=\inf\{ \lambda > 0 : \lambda u \in
\W_\1^s\}$ for $u \in V\cap\partial B_c(0)$. If we consider\footnote{By
$\partial[V\cap \partial B_c(0)]$ we mean the border of the $(d-1)$-dimensional
manifold $V\cap \partial B_c(0)$.}
$$
\eta:= \inf_{u \in \partial[V\cap\partial B_c(0)]} \phi(\bar{\lambda}_u u) >
\phi(\1)
$$ then the fact that $\phi(U(t))$ is strictly decreasing (see Proposition
\ref{Lyapunov}) allows us to shrink $V$ into a smaller neighborhood $V^*$ of
$c\1$ such that $\phi(v)=\eta$ for all $v \in \partial[V^*\cap \partial
B_c(0)]$. Let us also observe that since $\1$ is the only saddle point we can
take $V$ sufficiently small so as to guarantee that $\max\{ \phi(\lambda u) :
\lambda > 0\} \geq \eta$ for all $u \in \partial B_c(0)\setminus V^*$. Then if
we take the level curve $C_\eta = \{ x \in \R^d : \phi(x) = \eta \}$ every ray
$R_u$ with $u \in \partial B_c(0)\setminus V^*$ intersects $C_\eta$. With this
we may define for each $u \in \partial B_c(0)$
$$
\lambda_u^*=
\left\{\begin{array}{ll}
\bar{\lambda}_u & \mbox{ if }\,\,u \in V^*\\
\\
\inf\{ \lambda > 0 : \lambda u \in C_\eta\} & \mbox{ if }\,\,u \in B_c(0)\setminus V^*
\end{array}
\right..
$$ Notice that the application $u \mapsto \lambda^*_u$ is continuous. Due to
this fact, if $\tilde{G}:=\{ \lambda u : 0 \leq \lambda < \lambda^*_u\,,\, u \in
\partial B_c(0)\}$ then $\partial\tilde{G} =\{ \lambda^*_u u : u \in \partial
B_c(0)\}$. To finish the construction of our domain we must make a slight radial
expansion of $\tilde{G}$, i.e., for $\alpha > 0$ consider $G$ defined by the
formula
$$
G:=\{\lambda u : 0 \leq \lambda < (1+\alpha)\lambda^*_u\,,\, u \in \partial B_c(0)\}.
$$
Let us observe that Theorem \ref{divisionposta} insures that $G$ verifies
condition $(1)$. Since $\lambda^*_u > 1$ for all $u \in \partial B_c(0)$ then it
must also verify $(2)$. Also, if we define $\partial^1:=\{(1+\alpha)\lambda^*(u)
: \lambda^*(u) u \in \overline{V^*}\}$ then $\partial^1$ is closed and if
$\alpha > 0$ is taken small enough then $(3)$ holds. Finally, due to Theorem
\ref{divisionposta} we have $\partial^1 \subset D_e$ and so $(4)$ is verified. 
See Figure \ref{Gaux}.

\psfrag{0}{$0$}
\psfrag{1}{\vspace{-45pt}$\1$}
\psfrag{W1s}{$\W_\1^s$}
\psfrag{phieta}{$C_\eta$}

\subsection{The escape from $G$}

The behavior of the explosion time for initial data $u\in D_0$ is proved by
showing that, with overwhelming probability as $\varepsilon \to 0$, the
stochastic system describes the following path:
\begin{enumerate}
\item [(i)] It enters a neighborhood of the origin $B_c(0)$ in before a finite
time $T$ that does not depend on $\ve$.
\item [(ii)] Once in $B_c(0)$ the system remains in $G$ for a time of order
$e^{{\Delta}/{\ve^2}}$ and then escapes from $G$ through $\partial^1$ since the
barrier imposed by the potential is the lowest there.
\item [(iii)] After escaping $G$ through $\partial^1$ the system explodes
before a finite time $\tau$ which does not depend on $\ve$.
\end{enumerate}

The fact that the domain $G$ is bounded allows us to assume that $b$ is
globally Lipschitz if we wish to study the behavior of our system while it
remains inside $G$. Indeed, we may take $n_0 \in \N$ such that $G \subset
B_{n_0}(0)$ and study the behavior of the solution to \eqref{eqtruncada} since
it coincides with our process until the escape from $G$. Then we can proceed as
in the double-well potential case to obtain the following results (see \cite[pp
295--300]{OV} for their proofs). Hereafter, $B_c(0)$ denotes the neighborhood of
the origin highlighted in the construction of $G$ in the previous section.

\begin{teo}\label{ecotsuplema1} Given $\delta > 0$ we have
\begin{equation*}
\lim_{\varepsilon \rightarrow \: 0 }\:\: \sup_{u \in \overline B_c(0)} P_{u}
\Big( e^{\frac{\Delta - \delta}{\varepsilon^{2}}} < \tau_{\varepsilon}(\partial
G) < e^{\frac{\Delta + \delta}{\varepsilon^{2}}}\Big) = 1.
\end{equation*}
\end{teo}

\begin{teo}\label{teoescapedeg} The stochastic system verifies
$$
\lim_{\varepsilon \rightarrow 0} \sup_{u \in \overline{B}_c(0)} P_u \big(U^{\varepsilon}(\tau_\varepsilon(\partial G)) \notin \partial^1\big) = 0.
$$
\end{teo}

From these two theorems we can obtain the following useful corollary.

\begin{cor}\label{ecotsupcor1} For any $\delta > 0$ we have
\begin{equation*}
\lim_{\varepsilon \rightarrow 0} \sup_{u \in \overline{B}_{c}(0)} P_u
\Big(\tau_{\varepsilon}(\partial^1 ) > e^{\frac{\Delta +
\delta}{\varepsilon^2}}\Big)=0.
\end{equation*}
\end{cor}

\begin{proof} One can easily check that
$$
\sup_{u \in \overline{B}_c(0)} P_u \Big(\tau_\varepsilon(\partial^1) > e^{\frac{\Delta
+ \delta}{\varepsilon^2}}\Big) \leq \sup_{u \in \overline{B}_c(0)} P_u \Big(\tau_{\varepsilon}(\partial G) \geq e^{\frac{\Delta
+ \delta}{\varepsilon^2}}\Big) \:+ \sup_{u \in \overline{B}_c(0)}P_u \Big(U^\varepsilon_{\tau_\varepsilon(\partial G)} \notin \partial^1\Big).
$$
\end{proof}

Concerning the asymptotic distribution of $\tau_\ve(\partial G)$ we can obtain the following result.

\begin{teo}\label{nescapeteo1} Let $\gamma_\varepsilon
> 0$ be defined by the relation
\begin{equation*}
P_0 ( \tau_{\varepsilon} (\partial G) > \gamma_\varepsilon ) = e^{-1}.
\end{equation*}
Then there exists $\rho > 0$ such that for all $t\geq 0$ we have
\begin{equation*}
\lim_{\varepsilon \rightarrow 0} \sup_{u \in
\overline{B}_{\rho}(0)} |P_u (\tau_\varepsilon (\partial G) >
t\gamma_\varepsilon ) - e^{-t}| = 0.
\end{equation*}
\end{teo}

\subsection{Bounds for the explosion time}

This section is devoted to the lower and upper bounds for the explosion time.
More precisely, in this section we show that given $\delta > 0$, for all $u \in
D_0$ one has
\begin{equation*}\label{ncotainferior1}
\lim_{\varepsilon \rightarrow 0} P_u \Big (\tau_\varepsilon < e^{\frac{\Delta - \delta}{\varepsilon^2}}\Big) = 0
\end{equation*}and
\begin{equation*}\label{ncotasuperior1}
\lim_{\varepsilon \rightarrow 0} P_u \Big (\tau_\varepsilon > e^{\frac{\Delta + \delta}{\varepsilon^2}}\Big) = 0,
\end{equation*} where the convergence can be taken uniform over compact subsets
of $D_0$. The proofs of these bounds essentially follow \cite{OV}, where
analogous bounds are given for the tunneling time. However, unlike the
double-well potential model, the use of localization techniques becomes
necessary at some points throughout our work. We begin first with the lower
bound.

\begin{prop}\label{nteocotainferior0} Given $\delta > 0$ and $u \in D_0$ we have
\begin{equation}\label{cotainf001}
\lim_{\varepsilon \rightarrow 0} P_u \Big (\tau_\varepsilon < e^{\frac{\Delta - \delta}{\varepsilon^2}}\Big) = 0.
\end{equation}
Moreover, the convergence is uniform over compact subsets of $D_0$.
\end{prop}

\begin{proof} First observe that since for $u\in G$ we have $P_u(\tau_\ve\ge
\tau_\ve(\partial G))=1$ then \eqref{cotainf001} holds uniformly over any small
neighborhood of the origin by Lemma \ref{ecotsuplema1}. Next, we generalize the
result for any $u \in D_0$. For each $u \in D_0$ there exist $T_u > 0$,
$\delta_u > 0$ and $n_u \in \N$ such that the deterministic system beginning at
$u$ reaches $B_{\frac{\rho}{2}}(0)$ before $T_u$, remaining in $B_{n_u}(0)$ and
at a distance $\delta_u$ from $\partial B_{n_u}(0)$ on $[0, T_u]$.
It follows that $U^{n_u,\,u}$ does so as well. From this we obtain
\begin{align*}
P_u\big(\tau_{\varepsilon}(\overline{B}_{\rho}(0)) > T_u\big) &\leq P_u\big(
\min \{\tau_{\varepsilon}^{n_u}\,,\,\tau_{\varepsilon}(\overline{B}_{\rho}(0))\}
> T_u \big) + P_u\big(\tau_{\varepsilon}^{n_u} \leq T_u \big)\\
\\
& \leq P_u \Big(\sup_{0\le t \le T_u}|U^{n_u,\,\varepsilon}(t)- U^{n_u}(t)| > \frac{\rho}{2} \,\Big) + P_u \Big(\sup_{0\le t \le T_u} |U^{n_u,\,\varepsilon}(t)- U^{n_u}(t)| > \frac{\,\,\delta_u}{2} \Big).
\end{align*}
Using estimation (\ref{grandes1}) for the family
$\big( U^{n_u,\,u,\,\varepsilon}\big)_{\varepsilon > 0}$
we conclude
\begin{equation}
\label{Tu.menor.taueps}
\lim_{\varepsilon \rightarrow 0} P_u\big(\tau_{\varepsilon}(\overline{B}_{\rho}(0)) > T_u\big) = 0.
\end{equation}
Therefore, if we write
$$
P_u \Big( \tau_\varepsilon < e^{\frac{\Delta - \delta}{\varepsilon^2}} \Big)
\leq P_u \Big( \tau_{\varepsilon}(\overline{B}_{\rho}(0)) < \tau_\varepsilon <
e^{\frac{\Delta - \delta}{\varepsilon^2}} \Big) + P_u (\tau_\varepsilon \leq T_u
) + P_u (  \tau_{\varepsilon}(\overline{B}_{\rho}(0)) > T_u ),
$$
then the last two terms on the right tend to zero when $\varepsilon \rightarrow
0$ as a consequence of what we stated above. By the strong Markov property for
$U^{u,\,\varepsilon}$ we have
$$
P_u \Big( \tau_{\varepsilon}(\overline{B}_{\rho}(0)) < \tau_\varepsilon < e^{\frac{\Delta - \delta}{\varepsilon^2}} \Big) \leq \sup_{y  \in \overline{B}_\rho(0)} P_y \Big( \tau_\varepsilon < e^{\frac{\Delta - \delta}{\varepsilon^2}} \Big) \leq \sup_{y  \in \overline{B}_\rho(0)} P_y \Big( \tau^{n_0}(\partial G) < e^{\frac{\Delta - \delta}{\varepsilon^2}} \Big)
$$where $n_0$ is taken as in the first step. Since the rightmost term tends to zero by Lemma \ref{ecotsuplema1} we conclude the result for arbitrary $u \in D_0$. The uniform convergence over compact subsets $\mathcal{K}$ of $D_0$ is proved in a similar fashion by taking $\delta_u$ and $T_u$ uniformly over $\mathcal{K}$ as in Proposition \ref{convsup}.
\end{proof}

Now we turn to the proof of the upper bound. As we stated before, when studying
the behavior of the stochastic system under initial conditions $u \in G$ and for
small $\varepsilon > 0$ we typically observe that the process $U^{u,
\,\varepsilon}$ escapes from $G$ through $\partial^1$ since the cost imposed by
the potential is the lowest there. Once in $\partial^1$ the influence of noise
becomes negligible and the process then describes a path similar to the
deterministic trajectory until exploding in a finite time. We formalize this
statement in the following proposition.

\begin{prop}\label{teohard}
There exists $T_0 > 0$ such that
\begin{equation*}\label{explosionfrontera}
\lim_{\varepsilon \rightarrow 0} \sup_{u \in \partial^1} P_u ( \tau_\varepsilon > T_0 ) = 0.
\end{equation*}
\end{prop}

\begin{proof}
Since $\partial^1$ is a compact set contained in $D_e$, the proof follows from
Proposition \ref{convsup} and the fact that  $\sup_{u \in \partial^1} \tau_u^0 <
+\infty$.
\end{proof}

With this proposition we are able to conclude the upper bound. 

\begin{prop} For each $\delta > 0$ and $u \in D_0$ we have
\begin{equation*}
\lim_{\varepsilon \rightarrow 0} P_{u} \Big( \tau_{\varepsilon} >
e^{\frac{\Delta + \delta}{\varepsilon^{2}}}\Big) = 0.
\end{equation*}
Moreover, the convergence is uniform over compact subsets of $D_0$.
\end{prop}

\begin{proof} We proceed in two steps.

\noindent \textbf{1}. We check that given $\delta > 0$ we get
\begin{equation}
\label{eteolim0}
\lim_{\varepsilon \rightarrow 0} \sup_{x \in \overline{B}_c(0)} P_x \Big( \tau_{\varepsilon} >
e^{\frac{\Delta + \delta}{\varepsilon^{2}}}\Big) = 0.
\end{equation}
It is not hard to show that for $\varepsilon > 0$ small enough the strong
Markov property yields
$$
\sup_{u \in \overline{B}_c(0)} P_u \Big(\tau_{\varepsilon} > e^{\frac{\Delta + \delta}{\varepsilon^{2}}} \Big) \leq \sup_{u \in \overline{B}_c(0)} P_u \Bigg(\tau_{\varepsilon} (\partial^1) >
e^{\frac{\Delta + \frac{\delta}{2}}{\varepsilon^{2}}}\Bigg) + \sup_{u \in \partial^1} P_u (\tau_{\varepsilon} > T_0 ) + \sup_{u\in \overline B_c(0)} P_u(U_{\tau_\ve(\partial G)}^\ve \notin \partial^1)
$$
where $T_0 > 0$ is taken as in Proposition \ref{teohard}. We finish this first step by observing that the right hand side converges to zero. Indeed, the first term does so by Corollary \ref{ecotsupcor1}, the second by Proposition \ref{teohard} and the third by Lemma \ref{teoescapedeg}.

\medskip
\noindent \textbf{2}. We now generalize the result for $u \in D_0$. This follows from the fact that
\begin{equation*}\label{eq16}
P_u \Big(\tau_{\varepsilon} > e^{\frac{\Delta + \delta}{\varepsilon^{2}}}\Big) \leq \sup_{u \in \overline{B}_{c}(0)} P_u \Bigg(\tau_{\varepsilon} > \:\frac{e^{\frac{\Delta + \delta}{\varepsilon^{2}}}}{2}\Bigg) + P_u \big(\tau_{\varepsilon}(\overline{B}_{c}(0)) > T_u \big)
\end{equation*}
by the strong Markov property. Observing that the first term on the right hand side of the equation tends to zero by \eqref{eteolim0} and that the second term does by \eqref{Tu.menor.taueps}, we obtain our result. The convergence over compact subsets of $D_0$ can be seen as in Proposition \ref{convsup}.
\end{proof}

\subsection{Asymptotic distribution of the explosion time}

\noindent Our main objective in this section is to prove the asymptotic memory loss of the normalized explosion time
$\frac{\tau_{\varepsilon}}{\beta_\varepsilon}$. The proof focuses on studying the escape from $G$. The asymptotic memory loss for $\tau_\ve$ can be deduced once we show that the time in which the process exits from $G$ and the explosion time are asymptotically similar. We formalize this last statement in the following proposition.

\begin{prop}\label{nescapelema3} There exists a positive constant $T_0$ such that for all $u \in D_0 \cap G$
\begin{equation*}
\lim_{\varepsilon \rightarrow 0} P_x ( \tau_\varepsilon > \tau_\varepsilon (\partial G) + T_0 ) = 0.
\end{equation*}
\end{prop}

\begin{proof} Let us observe that by the strong Markov property
$$
P_u ( \tau_\varepsilon > \tau_\varepsilon (\partial G) + T_0 ) \leq \sup_{y \in
\partial^1} P_y (\tau_{\varepsilon} > T_0 ) + P_u \big(\tau_\varepsilon
(\partial G) < \tau_\varepsilon (\overline{B}_c(0)) \big) + \sup_{u \in
\overline{B}_c(0)} P_u \big(U^{\varepsilon}_{\tau_\varepsilon(\partial G)}
\notin \partial^1\big).
$$ We can now conclude our desired result by the use of Proposition \ref{teohard} and Lemma \ref{teoescapedeg}.
\end{proof}

We are now ready to establish the asymptotic memory loss of the explosion time.
Having the former proposition at our disposal, the rest of the proof is very
similar to the one offered in the double-well potential model.  We emphasize
that the main difference with this case lies in how to show this last
proposition. In the double-well potential the corresponding statement to
Proposition \ref{nescapelema3} holds due to the fact that the tunneling time for
initial conditions in the deepest well is of order one. This can be easily
deduced from the Freidlin-Wentzell estimates. Analogously, in our model
Proposition \ref{nescapelema3} holds since now the explosion time for initial
data in $D_e$ is of order one. However, the lack of a global Lipschitz condition
forces us to proceed differently in order to show this last fact. We recall that
a proof of this is contained essentially in Proposition \ref{convsup}. We now
give a brief sketch of the rest of the proof of Theorem \ref{distribution} in
the following lines and refer to \cite{GOV} for further details.

\noindent\textbf{Sketch of proof of Theorem \ref{distribution}}.
\begin{enumerate}
\item We first check that, for $\rho > 0$ small enough, $\lim_{\varepsilon
\rightarrow 0} \sup_{u \in \overline{B}_{\rho}(0)} |P_u (\tau_\varepsilon
(\partial G) > t\beta_\varepsilon ) - e^{-t}| = 0.$ This is due to the fact that
$\lim_{\varepsilon \rightarrow 0} \frac{\beta_\varepsilon}{\gamma_\varepsilon} =
1$.
\item  Next, we prove that $P_0 (\tau_{\varepsilon} > t\beta_{\varepsilon}) =
e^{-t}$ for $t > 0.$ This is done with the help of Proposition
\ref{nescapelema3} and the previous step.
\item With the help of appropriate coupling techniques we establish the uniform
convergence over any small enough neighborhood of the origin.
\item Finally, by using the strong Markov property, we conclude the result for
arbitrary initial data $u \in D_0$.
\end{enumerate}

\begin{center}
 {\bf Acknowledgments}
\end{center}
We want to thank Inés Armendáriz, Daniel Carando and Julián Fernández Bonder
for interesting discussions that helped us build this article. 

Both authors are partially supported by Universidad de Buenos Aires under grant
X447, by ANPCyT PICT 2008-00315 and CONICET PIP 643.

\bibliographystyle{amsplain}
\bibliography{bibliografia}

\end{document}